\documentclass{amsart}%
\usepackage{amssymb,amsmath,amsfonts,latexsym,amsthm,geometry,graphicx}
\usepackage{amsmath}
\usepackage{amsfonts}
\usepackage{amssymb}
\usepackage{graphicx}%
\providecommand{\U}[1]{\protect\rule{.1in}{.1in}}
\providecommand{\U}[1]{\protect\rule{.1in}{.1in}}
\providecommand{\U}[1]{\protect\rule{.1in}{.1in}}
\newtheorem{theorem}{Theorem}[section]

\theoremstyle{definition}

\begin{document}
\title{On the upper bounds for the constants of the Hardy--Littlewood inequality }
\author[G. Ara\'{u}jo]{Gustavo Ara\'{u}jo}
\address{Departamento de Matem\'{a}tica \\
Universidade Federal da Para\'{\i}ba \\
58.051-900 - Jo\~{a}o Pessoa, Brazil.}
\email{gdasaraujo@gmail.com}
\author[D. Pellegrino]{Daniel Pellegrino}
\address{Departamento de Matem\'{a}tica \\
Universidade Federal da Para\'{\i}ba \\
58.051-900 - Jo\~{a}o Pessoa, Brazil.}
\email{pellegrino@pq.cnpq.br and dmpellegrino@gmail.com}
\author[D.D. Pereira da Silva]{Diogo Diniz P. da Silva e Silva}
\address{Unidade Academica de Matematica e Estatistica \\
Universidade Federal de Campina Grande \\
Caixa Postal 10044 \\
58.429-270 - Campina Grande- PB, Brazil.}
\email{diogo@dme.ufcg.edu.br}
\thanks{The authors are supported by CNPq Grant 313797/2013-7 - PVE - Linha 2}

\begin{abstract}
The best known upper estimates for the constants of the Hardy--Littlewood
inequality for $m$-linear forms on $\ell_{p}$ spaces are of the form $\left(
\sqrt{2}\right)  ^{m-1}.$ We present better estimates which depend on $p$ and
$m$. An interesting consequence is that if $p\geq m^{2}$ then the constants
have a subpolynomial growth as $m$ tends to infinity.

\end{abstract}
\maketitle

\section{Introduction}

Let $\mathbb{K}$ be $\mathbb{R}$ or $\mathbb{C}$. Given an integer $m\geq2$,
the Hardy--Littlewood inequality (see \cite{alb, hardy,pra}) asserts that for
$2m\leq p\leq\infty$ there exists a constant $C_{m,p}^{\mathbb{K}}\geq1$ such
that, for all continuous $m$--linear forms $T:\ell_{p}^{n}\times\cdots
\times\ell_{p}^{n}\rightarrow\mathbb{K}$ and all positive integers $n$,
\begin{equation}
\left(  \sum_{j_{1},...,j_{m}=1}^{n}\left\vert T(e_{j_{1}},...,e_{j_{m}%
})\right\vert ^{\frac{2mp}{mp+p-2m}}\right)  ^{\frac{mp+p-2m}{2mp}}\leq
C_{m,p}^{\mathbb{K}}\left\Vert T\right\Vert .\label{i99}%
\end{equation}
Using the generalized Kahane-Salem-Zygmund inequality (see \cite{alb}) one can
easily verify that the exponents $\frac{2mp}{mp+p-2m}$ are optimal. The case
$p=\infty$ recovers the classical Bohnenblust--Hille inequality (see
\cite{bh}). More precisely, it asserts that there exists a constant
$B_{\mathbb{K},m}^{\mathrm{mult}}$ such that for all continuous $m$--linear
forms $T:\ell_{\infty}^{n}\times\cdots\times\ell_{\infty}^{n}\rightarrow
\mathbb{K}$ and all positive integers $n$,
\begin{equation}
\left(  \sum_{j_{1},...,j_{m}=1}^{n}\left\vert T(e_{j_{1}},...,e_{j_{m}%
})\right\vert ^{\frac{2m}{m+1}}\right)  ^{\frac{m+1}{2m}}\leq B_{\mathbb{K}%
,m}^{\mathrm{mult}}\left\Vert T\right\Vert .\label{u88}%
\end{equation}
From \cite{bohr, 111} we know that $B_{\mathbb{K},m}^{\mathrm{mult}}$ has a
subpolynomial growth. On the other hand, the best known upper bounds for the
constants in (\ref{i99}) are $\left(  \sqrt{2}\right)  ^{m-1}$ (see \cite{alb,
n, dimant}). In this paper we show that $\left(  \sqrt{2}\right)  ^{m-1}$ can
be improved to%
\[
C_{m,p}^{\mathbb{R}}\leq\left(  \sqrt{2}\right)  ^{\frac{2m\left(  m-1\right)
}{p}}\left(  B_{\mathbb{R},m}^{\mathrm{mult}}\right)  ^{\frac{p-2m}{p}}%
\]
for real scalars and to%
\[
C_{m,p}^{\mathbb{C}}\leq\left(  \frac{2}{\sqrt{\pi}}\right)  ^{\frac
{2m(m-1)}{p}}\left(  B_{\mathbb{C},m}^{\mathrm{mult}}\right)  ^{\frac{p-2m}%
{p}}%
\]
for complex scalars. These estimates are quite better than $\left(  \sqrt
{2}\right)  ^{m-1}$ because $B_{\mathbb{K},m}^{\mathrm{mult}}$ has a
subpolynomial growth. Moreover, our estimates depend on $p$ and $m$ and catch
more subtle information. For instance, if $p\geq m^{2}$ we conclude that
$\left(  C_{m,p}^{\mathbb{K}}\right)  _{m=1}^{\infty}$ has a subpolynomial
growth. Our main result is the following:

\begin{theorem}
Let $m\geq2$ be a positive integer and $2m\leq p\leq\infty.$ Then, for all
continuous $m$--linear forms $T:\ell_{p}^{n}\times\cdots\times\ell_{p}%
^{n}\rightarrow\mathbb{K}$ and all positive integers $n$, we have
\[
\left(  \sum_{j_{1},...,j_{m}=1}^{n}\left\vert T(e_{j_{1}},...,e_{j_{m}%
})\right\vert ^{\frac{2mp}{mp+p-2m}}\right)  ^{\frac{mp+p-2m}{2mp}}\leq
C_{m,p}^{\mathbb{K}}\left\Vert T\right\Vert
\]
with%
\[
C_{m,p}^{\mathbb{R}}\leq\left(  \sqrt{2}\right)  ^{\frac{2m\left(  m-1\right)
}{p}}\left(  B_{\mathbb{R},m}^{\mathrm{mult}}\right)  ^{\frac{p-2m}{p}}%
\]
and%
\[
C_{m,p}^{\mathbb{C}}\leq\left(  \frac{2}{\sqrt{\pi}}\right)  ^{\frac{2m\left(
m-1\right)  }{p}}\left(  B_{\mathbb{C},m}^{\mathrm{mult}}\right)
^{\frac{p-2m}{p}}.
\]

\end{theorem}

\section{The proof}

We recall that the Khinchin inequality (see \cite{Di}) asserts that \ for any
$0<q<\infty$, there are positive constants $A_{q}$, $B_{q}$ such that
regardless of the scalar sequence $(a_{j})_{j=1}^{\infty}$ in $\ell_{2}$ we
have
\[
A_{q}\left(  \sum_{j=1}^{\infty}|a_{j}|^{2}\right)  ^{\frac{1}{2}}\leq\left(
\int_{0}^{1}\left\vert \sum_{j=1}^{\infty}a_{j}r_{j}(t)\right\vert
^{q}dt\right)  ^{\frac{1}{q}}\leq B_{q}\left(  \sum_{j=1}^{\infty}|a_{j}%
|^{2}\right)  ^{\frac{1}{2}},
\]
where $r_{j}$ are the Rademacher functions. More generally, from the above
inequality together with the Minkowski inequality we know that
\begin{equation}
A_{q}^{m}\left(  \sum_{j_{1}...j_{m}=1}^{\infty}|a_{j_{1}...j_{m}}%
|^{2}\right)  ^{\frac{1}{2}}\leq\left(  \int_{I}\left\vert \sum_{j=1}^{\infty
}a_{j_{1}...j_{m}}r_{j_{1}}(t_{1})...r_{j_{m}}(t_{m})\right\vert ^{q}%
dt_{1}...dt_{m}\right)  ^{\frac{1}{q}}\leq B_{q}^{m}\left(  \sum
_{j_{1}...j_{m}=1}^{\infty}|a_{j_{1}...j_{m}}|^{2}\right)  ^{\frac{1}{2}%
}\label{9090}%
\end{equation}
for $I=[0,1]^{m}$ and all $\left(  a_{j_{1}....j_{m}}\right)  _{j_{1}%
,...,j_{m}=1}^{\infty}$ in $\ell_{2}$. \bigskip The notation of the constant
$A_{q}$ above will be used in all this paper.

Let $1\leq s\leq2$ and
\[
\lambda_{0}=\frac{2s}{ms+s-2m+2}.
\]
Since%
\[
\frac{m-1}{s}+\frac{1}{\lambda_{0}}=\frac{m+1}{2},
\]
from the generalized Bohnenblust--Hille inequality (see \cite{alb}) we know
that there is a constant $C_{m}\geq1$ such that for all $m$-linear forms
$T:\ell_{\infty}^{n}\times\cdots\times\ell_{\infty}^{n}\rightarrow\mathbb{K}$
we have%

\begin{equation}
\left(  \sum\limits_{j_{i}=1}^{n}\left(  \sum\limits_{\widehat{j_{i}}=1}%
^{n}\left\vert T\left(  e_{j_{1}},...,e_{j_{m}}\right)  \right\vert
^{s}\right)  ^{\frac{1}{s}\lambda_{0}}\right)  ^{\frac{1}{\lambda_{0}}}\leq
C_{m}\left\Vert T\right\Vert . \label{78}%
\end{equation}
Above, $\sum\limits_{\widehat{j_{i}}=1}^{n}$ means the sum over all $j_{k}$
for all $k\neq i.$ If we choose%
\[%
\begin{array}
[c]{ll}%
\displaystyle s=\frac{2mp}{mp+p-2m} & \text{ if }p<\infty,\vspace{0.2cm}\\
\displaystyle s=\frac{2m}{m+1} & \text{ if }p=\infty,
\end{array}
\]
we have%
\[
\lambda_{0}\leq s\leq2,
\]
and $\lambda_{0}=s$ when $p=\infty.$

The multiple exponent%
\[
\left(  \lambda_{0},s,s,...,s\right)
\]
can be obtained by interpolating the multiple exponents $\left(
1,2...,2\right)  \vspace{0.2cm}$ and $\left(  \frac{2m}{m+1},...,\frac
{2m}{m+1}\right)  $ with, respectively,
\[%
\begin{array}
[c]{c}%
\displaystyle\theta_{1}=2\left(  \frac{1}{\lambda_{0}}-\frac{1}{s}\right)
\vspace{0.2cm}\\
\displaystyle\theta_{2}=m\left(  \frac{2}{s}-1\right)  ,
\end{array}
\]
in the sense of \cite{alb}.

It is thus important to control the constants associated to the multiple
exponents $\left(  1,2...,2\right)  \vspace{0.2cm}$ and $\left(  \frac
{2m}{m+1},...,\frac{2m}{m+1}\right)  .$ The exponent $\left(  \frac{2m}%
{m+1},...,\frac{2m}{m+1}\right)  $ is the classical exponent of the
Bohnenblust--Hille inequality and the estimate of the constant associated to
$\left(  1,2...,2\right)  $ is well-known (we present the details for the sake
of completeness). In fact, in general, for the exponent $\left(  \frac
{2k}{k+1},...,\frac{2k}{k+1},2,...,2\right)  $ (with $\frac{2k}{k+1}$ repeated
$k$ times and $2$ repeated $m-k$ times), using the multiple Khinchin
inequality (\ref{9090}) , we have, for all $m$-linear forms $T:\ell_{\infty
}^{n}\times\cdots\times\ell_{\infty}^{n}\rightarrow\mathbb{K}$,
\begin{align*}
&  \left(  \sum\limits_{j_{1},...,j_{k}=1}^{n}\left(  \sum\limits_{j_{k+1}%
,...,j_{m}=1}^{n}\left\vert T\left(  e_{j_{1}},...,e_{j_{m}}\right)
\right\vert ^{2}\right)  ^{\frac{1}{2}\frac{2k}{k+1}}\right)  ^{\frac{k+1}%
{2k}}\\
&  \leq\left(  \sum\limits_{j_{1},...,j_{k}=1}^{n}\left(  A_{\frac{2k}{k+1}%
}^{-\left(  m-k\right)  }\left(  \int_{[0,1]^{m-k}}\left\vert \sum
\limits_{j_{k+1},...,j_{m}=1}^{n}r_{j_{k+1}}(t_{k+1})...r_{j_{m}}%
(t_{m})T\left(  e_{j_{1}},...,e_{j_{m}}\right)  \right\vert ^{\frac{2k}{k+1}%
}dt_{k+1}...dt_{m}\right)  ^{\frac{k+1}{2k}}\right)  ^{\frac{2k}{k+1}}\right)
^{\frac{k+1}{2k}}\\
&  =A_{\frac{2k}{k+1}}^{-\left(  m-k\right)  }\left(  \sum\limits_{j_{1}%
,...,j_{k}=1}^{n}\int_{[0,1]^{m-k}}\left\vert T\left(  e_{j_{1}},...,e_{j_{k}%
},\sum\limits_{j_{k+1}=1}^{n}r_{j_{k+1}}(t_{k+1})e_{j_{k+1}},...,\sum
\limits_{j_{m}=1}^{n}r_{j_{m}}(t_{m})e_{j_{m}}\right)  \right\vert ^{\frac
{2k}{k+1}}dt_{k+1}...dt_{m}\right)  ^{\frac{k+1}{2k}}\\
&  =A_{\frac{2k}{k+1}}^{-\left(  m-k\right)  }\left(  \int_{[0,1]^{m-k}}%
\sum\limits_{j_{1},...,j_{k}=1}^{n}\left\vert T\left(  e_{j_{1}},...,e_{j_{k}%
},\sum\limits_{j_{k+1}=1}^{n}r_{j_{k+1}}(t_{k+1})e_{j_{k+1}},...,\sum
\limits_{j_{m}=1}^{n}r_{j_{m}}(t_{m})e_{j_{m}}\right)  \right\vert ^{\frac
{2k}{k+1}}dt_{k+1}...dt_{m}\right)  ^{\frac{k+1}{2k}}\\
&  \leq A_{\frac{2k}{k+1}}^{-\left(  m-k\right)  }\sup_{t_{k+1},...,t_{m}%
\in\lbrack0,1]}B_{\mathbb{K},k}^{\mathrm{mult}}\left\Vert T\left(
\ \cdot\ ,...,\ \cdot\ ,\sum\limits_{j_{k+1}=1}^{n}r_{j_{k+1}}(t_{k+1}%
)e_{j_{k+1}},...,\sum\limits_{j_{m}=1}^{n}r_{j_{m}}(t_{m})e_{j_{m}}\right)
\right\Vert \\
&  =A_{\frac{2k}{k+1}}^{-\left(  m-k\right)  }B_{\mathbb{K},k}^{\mathrm{mult}%
}\left\Vert T\right\Vert .
\end{align*}

So, choosing $k=1$, since $A_{1}=\left(  \sqrt{2}\right)  ^{-1}$ and
$B_{\mathbb{K},1}^{\mathrm{mult}}=1$ we conclude that the constant associated
to the multiple exponent $\left(  1,2,...,2\right)  $ is $\left(  \sqrt
{2}\right)  ^{m-1}.$

Therefore, the optimal constant associated to the multiple exponent
\[
\left(  \lambda_{0},s,s,...,s\right)
\]
is less or equal (for real scalars) than%
\[
\left(  \left(  \sqrt{2}\right)  ^{m-1}\right)  ^{2\left(  \frac{1}%
{\lambda_{0}}-\frac{1}{s}\right)  }\left(  B_{\mathbb{R},m}^{\mathrm{mult}%
}\right)  ^{m\left(  \frac{2}{s}-1\right)  }%
\]
i.e.,%
\begin{equation}
C_{m}\leq\left(  \sqrt{2}\right)  ^{\frac{2m\left(  m-1\right)  }{p}}\left(
B_{\mathbb{R},m}^{\mathrm{mult}}\right)  ^{\frac{p-2m}{p}}.\label{9898}%
\end{equation}
More precisely, (\ref{78}) is valid with $C_{m}$ as above. For complex scalars
we can use the Khinchin inequality for Steinhaus variables and replace
$\sqrt{2}$ by $\frac{2}{\sqrt{\pi}}$ as in \cite{ddss}.

Let%
\[
\lambda_{j}=\frac{\lambda_{0}p}{p-\lambda_{0}j}%
\]
for all $j=0,....,m.$ Note that%
\[
\lambda_{m}=s
\]
and that
\[
\left(  \frac{p}{\lambda_{j}}\right)  ^{\ast}=\frac{\lambda_{j+1}}{\lambda
_{j}}%
\]
for all $j=0,...,m-1$.

Let us suppose that $1\leq k\leq m$ and that
\[
\left(  \sum_{j_{i}=1}^{n}\left(  \sum_{\widehat{j_{i}}=1}^{n}\left\vert
T(e_{j_{1}},...,e_{j_{m}})\right\vert ^{s}\right)  ^{\frac{1}{s}\lambda_{k-1}%
}\right)  ^{\frac{1}{\lambda_{k-1}}}\leq C_{m}\Vert T\Vert
\]
is true for all continuous $m$--linear forms $T:\underbrace{\ell_{p}^{n}%
\times\cdots\times\ell_{p}^{n}}_{k-1\ \mathrm{times}}\times\ell_{\infty}%
^{n}\times\cdots\times\ell_{\infty}^{n}\rightarrow\mathbb{K}$ and for all
$i=1,...,m.$ Let us prove that
\[
\left(  \sum_{j_{i}=1}^{n}\left(  \sum_{\widehat{j_{i}}=1}^{n}\left\vert
T(e_{j_{1}},...,e_{j_{m}})\right\vert ^{s}\right)  ^{\frac{1}{s}\lambda_{k}%
}\right)  ^{\frac{1}{\lambda_{k}}}\leq C_{m}\Vert T\Vert
\]
for all continuous $m$--linear forms $T:\underbrace{\ell_{p}^{n}\times
\cdots\times\ell_{p}^{n}}_{k\ \mathrm{times}}\times\ell_{\infty}^{n}%
\times\cdots\times\ell_{\infty}^{n}\rightarrow\mathbb{K}$ and for all
$i=1,...,m$.

The initial case (the case $k=0$) is precisely (\ref{78}) with $C_{m}$ as in
(\ref{9898}).

Consider
\[
T\in\mathcal{L}(\underbrace{\ell_{p}^{n},...,\ell_{p}^{n}}_{k\ \mathrm{times}%
},\ell_{\infty}^{n},...,\ell_{\infty}^{n};\mathbb{K})
\]
and for each $x\in B_{\ell_{p}^{n}}$ define
\[%
\begin{array}
[c]{ccccl}%
T^{(x)} & : & \underbrace{\ell_{p}^{n}\times\cdots\times\ell_{p}^{n}%
}_{k-1\ \mathrm{times}}\times\ell_{\infty}^{n}\times\cdots\times\ell_{\infty
}^{n} & \rightarrow & \mathbb{K}\\
&  & (z^{(1)},...,z^{(m)}) & \mapsto & T(z^{(1)},...,z^{(k-1)},xz^{(k)}%
,z^{(k+1)},...,z^{(m)}),
\end{array}
\]
with $xz^{(k)}=(x_{j}z_{j}^{(k)})_{j=1}^{n}$. Observe that
\[
\Vert T\Vert=\sup\{\Vert T^{(x)}\Vert:x\in B_{\ell_{p}^{n}}\}.
\]
By applying the induction hypothesis to $T^{(x)}$, we obtain
\begin{equation}%
\begin{array}
[c]{l}%
\displaystyle\left(  \sum_{j_{i}=1}^{n}\left(  \sum_{\widehat{j_{i}}=1}%
^{n}\left\vert T\left(  e_{j_{1}},...,e_{j_{m}}\right)  \right\vert
^{s}\left\vert x_{j_{k}}\right\vert ^{s}\right)  ^{\frac{1}{s}\lambda_{k-1}%
}\right)  ^{\frac{1}{\lambda_{k-1}}}\vspace{0.2cm}\\
=\displaystyle\left(  \sum_{j_{i}=1}^{n}\left(  \sum_{\widehat{j_{i}}=1}%
^{n}\left\vert T\left(  e_{j_{1}},...,e_{j_{k-1}},xe_{j_{k}},e_{j_{k+1}%
},...,e_{j_{m}}\right)  \right\vert ^{s}\right)  ^{\frac{1}{s}\lambda_{k-1}%
}\right)  ^{\frac{1}{\lambda_{k-1}}}\vspace{0.2cm}\\
=\displaystyle\left(  \sum_{j_{i}=1}^{n}\left(  \sum_{\widehat{j_{i}}=1}%
^{n}\left\vert T^{(x)}\left(  e_{j_{1}},...,e_{j_{m}}\right)  \right\vert
^{s}\right)  ^{\frac{1}{s}\lambda_{k-1}}\right)  ^{\frac{1}{\lambda_{k-1}}%
}\vspace{0.2cm}\\
\leq C_{m}\Vert T^{(x)}\Vert\vspace{0.2cm}\\
\leq C_{m}\Vert T\Vert
\end{array}
\label{guga030}%
\end{equation}
for all $i=1,...,m.$

We will analyze two cases:

\bigskip

\noindent1) $i=k$

\bigskip

Since
\[
\left(  \frac{p}{\lambda_{j-1}}\right)  ^{\ast}=\frac{\lambda_{j}}%
{\lambda_{j-1}}%
\]
for all $j=1,...,m$, we conclude that
\[%
\begin{array}
[c]{l}%
\displaystyle\left(  \sum_{j_{k}=1}^{n}\left(  \sum_{\widehat{j_{k}}=1}%
^{n}\left\vert T\left(  e_{j_{1}},...,e_{j_{m}}\right)  \right\vert
^{s}\right)  ^{\frac{1}{s}\lambda_{k}}\right)  ^{\frac{1}{\lambda_{k}}}%
\vspace{0.2cm}\\
\displaystyle=\displaystyle\left(  \sum_{j_{k}=1}^{n}\left(  \sum
_{\widehat{j_{k}}=1}^{n}\left\vert T\left(  e_{j_{1}},...,e_{j_{m}}\right)
\right\vert ^{s}\right)  ^{\frac{1}{s}\lambda_{k-1}\left(  \frac{p}%
{\lambda_{k-1}}\right)  ^{\ast}}\right)  ^{\frac{1}{\lambda_{k-1}}\frac
{1}{\left(  \frac{p}{\lambda_{k-1}}\right)  ^{\ast}}} \vspace{0.2cm}\\
\displaystyle=\left(  \left\Vert \left(  \left(  \sum_{\widehat{j_{k}}=1}%
^{n}\left\vert T\left(  e_{j_{1}},...,e_{j_{m}}\right)  \right\vert
^{s}\right)  ^{\frac{1}{s}\lambda_{k-1}}\right)  _{j_{k}=1}^{n}\right\Vert
_{\left(  \frac{p}{\lambda_{k-1}}\right)  ^{\ast}}\right)  ^{\frac{1}%
{\lambda_{k-1}}}\vspace{0.2cm}\\
\displaystyle=\left(  \sup_{y\in B_{\ell_{\frac{p}{\lambda_{k-1}}}^{n}}}%
\sum_{j_{k}=1}^{n}|y_{j_{k}}|\left(  \sum_{\widehat{j_{k}}=1}^{n}\left\vert
T\left(  e_{j_{1}},...,e_{j_{m}}\right)  \right\vert ^{s}\right)  ^{\frac
{1}{s}\lambda_{k-1}}\right)  ^{\frac{1}{\lambda_{k-1}}}\vspace{0.2cm}\\
\displaystyle=\left(  \sup_{x\in B_{\ell_{p}^{n}}}\sum_{j_{k}=1}^{n}|x_{j_{k}%
}|^{\lambda_{k-1}}\left(  \sum_{\widehat{j_{k}}=1}^{n}\left\vert T\left(
e_{j_{1}},...,e_{j_{m}}\right)  \right\vert ^{s}\right)  ^{\frac{1}{s}%
\lambda_{k-1}}\right)  ^{\frac{1}{\lambda_{k-1}}}\vspace{0.2cm}\\
\displaystyle=\sup_{x\in B_{\ell_{p}^{n}}}\left(  \sum_{j_{k}=1}^{n}\left(
\sum_{\widehat{j_{k}}=1}^{n}\left\vert T\left(  e_{j_{1}},...,e_{j_{m}%
}\right)  \right\vert ^{s}\left\vert x_{j_{k}}\right\vert ^{s}\right)
^{\frac{1}{s}\lambda_{k-1}}\right)  ^{\frac{1}{\lambda_{k-1}}}\vspace{0.2cm}\\
\displaystyle\leq C_{m}\Vert T\Vert.
\end{array}
\]
where the last inequality holds by \eqref{guga030}.

\bigskip

\noindent2) $i\neq k$

\bigskip

It is important to note that $\lambda_{k-1}<\lambda_{k}\leq s$. Denoting, for
$i=1,....,m,$
\[
S_{i}=\left(  \sum_{\widehat{j_{i}}=1}^{n}|T(e_{j_{1}},...,e_{j_{m}}%
)|^{s}\right)  ^{\frac{1}{s}}%
\]
we get
\[%
\begin{array}
[c]{l}%
\displaystyle\sum_{j_{i}=1}^{n}\left(  \sum_{\widehat{j_{i}}=1}^{n}%
|T(e_{j_{1}},...,e_{j_{m}})|^{s}\right)  ^{\frac{1}{s}\lambda_{k}}=\sum
_{j_{i}=1}^{n}S_{i}^{\lambda_{k}}=\sum_{j_{i}=1}^{n}S_{i}^{\lambda_{k}-s}%
S_{i}^{s}\vspace{0.2cm}\\
\displaystyle=\sum_{j_{i}=1}^{n}\sum_{\widehat{j_{i}}=1}^{n}\frac{|T(e_{j_{1}%
},...,e_{j_{m}})|^{s}}{S_{i}^{s-\lambda_{k}}}=\sum_{j_{k}=1}^{n}\sum
_{\widehat{j_{k}}=1}^{n}\frac{|T(e_{j_{1}},...,e_{j_{m}})|^{s}}{S_{i}%
^{s-\lambda_{k}}}\vspace{0.2cm}\\
\displaystyle=\sum_{j_{k}=1}^{n}\sum_{\widehat{j_{k}}=1}^{n}\frac{|T(e_{j_{1}%
},...,e_{j_{m}})|^{\frac{s(s-\lambda_{k})}{s-\lambda_{k-1}}}}{S_{i}%
^{s-\lambda_{k}}}|T(e_{j_{1}},...,e_{j_{m}})|^{\frac{s(\lambda_{k}%
-\lambda_{k-1})}{s-\lambda_{k-1}}}.
\end{array}
\]
Therefore, using H\"{o}lder's inequality twice we obtain
\begin{equation}%
\begin{array}
[c]{l}%
\displaystyle\sum_{j_{i}=1}^{n}\left(  \sum_{\widehat{j_{i}}=1}^{n}%
|T(e_{j_{1}},...,e_{j_{m}})|^{s}\right)  ^{\frac{1}{s}\lambda_{k}}
\vspace{0.2cm}\\
\displaystyle\leq\sum_{j_{k}=1}^{n}\left(  \sum_{\widehat{j_{k}}=1}^{n}%
\frac{|T(e_{j_{1}},...,e_{j_{m}})|^{s}}{S_{i}^{s-\lambda_{k-1}}}\right)
^{\frac{s-\lambda_{k}}{s-\lambda_{k-1}}}\left(  \sum_{\widehat{j_{k}}=1}%
^{n}|T(e_{j_{1}},...,e_{j_{m}})|^{s}\right)  ^{\frac{\lambda_{k}-\lambda
_{k-1}}{s-\lambda_{k-1}}}\vspace{0.2cm}\\
\displaystyle\leq\left(  \sum_{j_{k}=1}^{n}\left(  \sum_{\widehat{j_{k}}%
=1}^{n}\frac{|T(e_{j_{1}},...,e_{j_{m}})|^{s}}{S_{i}^{s-\lambda_{k-1}}%
}\right)  ^{\frac{\lambda_{k}}{\lambda_{k-1}}}\right)  ^{\frac{\lambda_{k-1}%
}{\lambda_{k}}\cdot\frac{s-\lambda_{k}}{s-\lambda_{k-1}}}\left(  \sum
_{j_{k}=1}^{n}\left(  \sum_{\widehat{j_{k}}=1}^{n}|T(e_{j_{1}},...,e_{j_{m}%
})|^{s}\right)  ^{\frac{1}{s}\lambda_{k}}\right)  ^{\frac{1}{\lambda_{k}}%
\cdot\frac{(\lambda_{k}-\lambda_{k-1})s}{s-\lambda_{k-1}}}.
\end{array}
\label{huhu}%
\end{equation}
We know from the case $i=k$ that
\begin{equation}
\displaystyle\left(  \sum_{j_{k}=1}^{n}\left(  \sum_{\widehat{j_{k}}=1}%
^{n}|T(e_{j_{1}},...,e_{j_{m}})|^{s}\right)  ^{\frac{1}{s}\lambda_{k}}\right)
^{\frac{1}{\lambda_{k}}\cdot\frac{(\lambda_{k}-\lambda_{k-1})s}{s-\lambda
_{k-1}}}\leq\left(  C_{m}\Vert T\Vert\right)  ^{\frac{(\lambda_{k}%
-\lambda_{k-1})s}{s-\lambda_{k-1}}}. \label{huhi}%
\end{equation}
Now we investigate the first factor in (\ref{huhu}). From H\"{o}lder's inequality and
\eqref{guga030} it follows that
\begin{equation}%
\begin{array}
[c]{l}%
\displaystyle\left(  \sum_{j_{k}=1}^{n}\left(  \sum_{\widehat{j_{k}}=1}%
^{n}\frac{|T(e_{j_{1}},...,e_{j_{m}})|^{s}}{S_{i}^{s-\lambda_{k-1}}}\right)
^{\frac{\lambda_{k}}{\lambda_{k-1}}}\right)  ^{\frac{\lambda_{k-1}}%
{\lambda_{k}}}=\left\Vert \left(  \sum_{\widehat{j_{k}}}\frac{|T(e_{j_{1}%
},...,e_{j_{m}})|^{s}}{S_{i}^{s-\lambda_{k-1}}}\right)  _{j_{k}=1}%
^{n}\right\Vert _{\left(  \frac{p}{\lambda_{k-1}}\right)  ^{\ast}}%
\vspace{0.2cm}\\
\displaystyle=\sup_{y\in B_{\ell_{\frac{p}{\lambda_{k-1}}}^{n}}}\sum_{j_{k}%
=1}^{n}|y_{j_{k}}|\sum_{\widehat{j_{k}}=1}^{n}\frac{|T(e_{j_{1}},...,e_{j_{m}%
})|^{s}}{S_{i}^{s-\lambda_{k-1}}}=\sup_{x\in B_{\ell_{p}^{n}}}\sum_{j_{k}%
=1}^{n}\sum_{\widehat{j_{k}}=1}^{n}\frac{|T(e_{j_{1}},...,e_{j_{m}})|^{s}%
}{S_{i}^{s-\lambda_{k-1}}}|x_{j_{k}}|^{\lambda_{k-1}}\vspace{0.2cm}\\
\displaystyle=\sup_{x\in B_{\ell_{p}^{n}}}\sum_{j_{i}=1}^{n}\sum
_{\widehat{j_{i}}=1}^{n}\frac{|T(e_{j_{1}},...,e_{j_{m}})|^{s-\lambda_{k-1}}%
}{S_{i}^{s-\lambda_{k-1}}}|T(e_{j_{1}},...,e_{j_{m}})|^{\lambda_{k-1}%
}|x_{j_{k}}|^{\lambda_{k-1}}\vspace{0.2cm}\\
\displaystyle\leq\sup_{x\in B_{\ell_{p}^{n}}}\sum_{j_{i}=1}^{n}\left(
\sum_{\widehat{j_{i}}=1}^{n}\frac{|T(e_{j_{1}},...,e_{j_{m}})|^{s}}{S_{i}^{s}%
}\right)  ^{\frac{s-\lambda_{k-1}}{s}}\left(  \sum_{\widehat{j_{i}}=1}%
^{n}|T(e_{j_{1}},...,e_{j_{m}})|^{s}|x_{j_{k}}|^{s}\right)  ^{\frac{1}%
{s}\lambda_{k-1}}\vspace{0.2cm}\\
\displaystyle=\sup_{x\in B_{\ell_{p}^{n}}}\sum_{j_{i}=1}^{n}\left(
\sum_{\widehat{j_{i}}=1}^{n}|T(e_{j_{1}},...,e_{j_{m}})|^{s}|x_{j_{k}}%
|^{s}\right)  ^{\frac{1}{s}\lambda_{k-1}}\leq\left(  C_{m}\Vert T\Vert\right)
^{\lambda_{k-1}}.
\end{array}
\label{huho}%
\end{equation}
Replacing \eqref{huhi} and \eqref{huho} in \eqref{huhu} we finally conclude
that
\begin{align*}
\displaystyle\sum_{j_{i}=1}^{n}\left(  \sum_{\widehat{j_{i}}=1}^{n}%
|T(e_{j_{1}},...,e_{j_{m}})|^{s}\right)  ^{\frac{1}{s}\lambda_{k}}  &
\leq\left(  C_{m}\Vert T\Vert\right)  ^{\lambda_{k-1}\frac{s-\lambda_{k}%
}{s-\lambda_{k-1}}}\left(  C_{m}\Vert T\Vert\right)  ^{\frac{(\lambda
_{k}-\lambda_{k-1})s}{s-\lambda_{k-1}}}\\
&  =\left(  C_{m}\Vert T\Vert\right)  ^{\lambda_{k}}.
\end{align*}

Since $\lambda_{m}=s$ the proof is done.

\section{Constants with subpolynomial growth}

The optimal constants of the Khinchin's inequality (these constants are due to
U. Haagerup \cite{haage}) are
\[
A_{q}=\sqrt{2}\left(  \frac{\Gamma\left(  \frac{q+1}{2}\right)  }{\sqrt{\pi}%
}\right)  ^{\frac{1}{q}}%
\]
for $q>q_{0}\cong1.847$ and
\[
A_{q}=2^{\frac{1}{2}-\frac{1}{q}}%
\]
for $q\leq q_{0}$, where $q_{0}\in(0,2)$ is the unique real number satisfying
\[
\Gamma\left(  \frac{q_{0}+1}{2}\right)  =\frac{\sqrt{\pi}}{2}.
\]
For complex scalars if we use the Khinchin inequality for Steinhaus variables
we have
\[
A_{q}=\left(  \Gamma\left(  \frac{q+2}{2}\right)  \right)  ^{\frac{1}{q}}%
\]
for all $1\leq q<2$ (see \cite{konig}).

\bigskip The best known upper estimates for $B_{\mathbb{R},m}^{\mathrm{mult}}$
and $B_{\mathbb{C},m}^{\mathrm{mult}}$ (from \cite{bohr}) are%
\[
B_{\mathbb{K},m}^{\mathrm{mult}}\leq%
{\displaystyle\prod\limits_{j=2}^{m}}
A_{\frac{2j-2}{j}}^{-1}.
\]
Combining these results we have%
\[%
\begin{array}
[c]{c}%
C_{m,p}^{\mathbb{R}} \leq\left(  2^{\frac{4m^{2}-pm-2m}{2p-4m}+\frac
{446381}{55440}}%
{\displaystyle\prod\limits_{j=14}^{m}}
\left(  \frac{\Gamma\left(  \frac{3}{2}-\frac{1}{j}\right)  }{\sqrt{\pi}%
}\right)  ^{\frac{j}{2-2j}}\right)  ^{\frac{p-2m}{p}}\text{ for }m\geq14
\vspace{0.2cm}\\
C_{m,p}^{\mathbb{R}} \leq\left(  \sqrt{2}\right)  ^{\frac{2m\left(
m-1\right)  }{p}}\left(
{\displaystyle\prod\limits_{j=2}^{m}}
2^{\frac{1}{2j-2}}\right)  ^{\frac{p-2m}{p}}\text{ for }2\leq m\leq13
\end{array}
\]
and%
\[
C_{m,p}^{\mathbb{C}}\leq\left(  \frac{2}{\sqrt{\pi}}\right)  ^{\frac{2m\left(
m-1\right)  }{p}}\left(
{\displaystyle\prod\limits_{j=2}^{m}}
\Gamma\left(  2-\frac{1}{j}\right)  ^{\frac{j}{2-2j}}\right)  ^{\frac{p-2m}%
{p}}.
\]

From \cite{bohr} we know that there is a constant $\kappa>0$ such that%
\[%
\begin{array}
[c]{c}%
B_{\mathbb{R},m}^{\mathrm{mult}} \leq\kappa m^{\frac{2-\ln2-\gamma}{2}}<\kappa
m^{0.37}, \vspace{0.2cm}\\
B_{\mathbb{C},m}^{\mathrm{mult}} \leq\kappa m^{\frac{1-\gamma}{2}}<\kappa
m^{0.22},
\end{array}
\]
for all $m$, where $\gamma$ is the Euler-Mascheroni constant. We thus conclude
that if $p\geq m^{2}$ then $\left(  C_{m,p}^{\mathbb{K}}\right)
_{m=1}^{\infty}$ has a subpolynomial growth.

\end{document}